\documentclass[a4paper]{amsart}
\usepackage{amsmath,amssymb,stmaryrd,graphicx,hyperref,color,colortbl,array,fancyvrb}
\usepackage[all]{xy}
\DeclareFontShape{OT1}{cmtt}{bx}{n}{<5><6><7><8><9><10><10.95><12><14.4><17.28><20.74><24.88>cmttb10}{}
\definecolor{darkblue}{RGB}{0,0,128}
\hypersetup{colorlinks={true},linkcolor={darkblue},citecolor={darkblue},filecolor={darkblue},urlcolor={darkblue},pdfauthor={Lukas Lewark},pdftitle={sl3-foam homology calculations}}
\allowdisplaybreaks

\newtheorem{thm}{Theorem}[section] 
\newtheorem{lem}[thm]{Lemma}        
\newtheorem*{conj}{Conjecture}       
\theoremstyle{definition}
\newtheorem{defn}[thm]{Definition}

\newcommand{\holy}[1]{\ensuremath{#1}}
\newcommand{\foamho}{\texttt{FoamHo}}

\newcommand{\pa}[1]{\ensuremath{\mathcal{#1}}}
\newcommand{\category}[1]{\ensuremath{\mathsf{#1}}}

\DeclareMathOperator{\id}{id}
\DeclareMathOperator{\Hom}{Hom}

\newcommand{\qRing}{\ensuremath{\mathbb{Z}[q^{\pm 1}]}}
\newcommand{\mymatrixfour}[4]{\ensuremath{\Biggl(\begin{matrix}#1&#2\\[-2mm]#3&#4\end{matrix}\Biggr)}}
\newcommand{\mymatrixtwov}[2]{\ensuremath{\Biggl(\begin{matrix}#1\\[-2mm]#2\end{matrix}\Biggr)}}
\newcommand{\mymatrixtwoh}[2]{\ensuremath{\bigl(\begin{matrix}#1&#2\end{matrix}\bigr)}}
\newcommand{\mymatrixsimple}[1]{#1}

\newcommand{\leftright}[2]{\ensuremath{#1\,\cdot\,#2}}

\newlength{\lowerhalftmp}
\newcommand{\aFigpath}{figs}
\newcommand{\lowerhalfx}[1]{\settoheight{\lowerhalftmp}{#1}\addtolength{\lowerhalftmp}{-1.2ex}\raisebox{-.5\lowerhalftmp}{#1}}

\newcommand{\avcfig}[2][]{\lowerhalfx{\includegraphics[#1]{\aFigpath/#2.pdf}}}
\newcommand{\afig}[2][]{{\includegraphics[#1]{\aFigpath/#2.pdf}}}
\newcommand{\myhfill}{\rule{0pt}{0pt}\hfill\rule{0pt}{0pt}}
\newcommand{\myparagraph}[1]{\paragraph{\bf #1}\medskip\quad}
\begin{document}
\title{$\mathfrak{sl}_3$-foam homology calculations}
\author{Lukas Lewark}
\date{\today}
\begin{abstract}
We exhibit a certain infinite family of three-stranded
quasi-alternating pretzel knots which are counterexamples to Lobb's conjecture
that the $\mathfrak{sl}_3$--knot concordance invariant $s_3$ (suitably normalised)
should be equal to the Rasmussen invariant $s_2$.
For this family, $|s_3| < |s_2|$.
However, we also find other knots for which $|s_3| > |s_2|$.
The main tool is an implementation of Morrison and Nieh's algorithm to calculate
Khovanov's $\mathfrak{sl}_3$--foam link homology. Our C++-program is fast enough
to calculate the integral homology of e.g. the $(6,5)$--torus knot in six minutes.
Furthermore, we propose a potential improvement of the algorithm
by gluing sub-tangles in a more flexible way.
%

\end{abstract}
\maketitle
\tableofcontents
\bibliographystyle{myamsalpha}
\section{Introduction}
During the last decade, $\mathfrak{sl}_n$--polynomials were categorified
one after the other:
beginning, of course, with Khovanov's categorification of the Jones polynomial \cite{khovanov}
using cobordisms, followed by his categorification of the $\mathfrak{sl}_3$--polynomial \cite{khovanovsl3}
using foams; and finally, Khovanov and Rozansky's categorification of the $\mathfrak{sl}_N$--polynomials
for arbitrary $N$
\cite{roz} and of the Homflypt--polynomial itself \cite{roz2} (see also Khovanov \cite{khsoergel} and Webster \cite{homflycan})
using matrix factorisations. 

All these homologies are completely combinatorial in nature -- unlike e.g. the original knot Floer homology --
meaning that their definition is in itself a description how to compute them.
By hand, this direct way of computation is hardly feasible for any but the smallest knots,
and using skein long exact sequences and criteria for thinness is much more efficient (see e.g.
a computation of Mackaay and Vaz \cite{terasaka}).
On a computer, however, even the straight-forward method, as implemented with some tweaks
in the program \texttt{KhoHo} \cite{khoho} by Shumakovitch, could already compute
the $\mathfrak{sl}_2$--homology of knots with up to ca. $20$ crossings.
Matrix factorisations, on the other hand, are more reluctant to efficient treatment by
a computer. Nevertheless, Carqueville and Murfet \cite{carqueville}
wrote a program that is able to compute $\mathfrak{sl}_N$--homology of links with
up to six crossings and $N \approx 18$.
Webster's program \cite{krm2} computes Homflypt-homology of short braids,
but is largely limited to 3--stranded ones.

Yet it is desirable to be able to compute the Khovanov--Rozansky homologies of much larger knots, since some phenomena
require a certain complexity of the knot to occur; see e.g. Hedden and Ording's 15--crossing knot
which demonstrates that the Rasmussen and Floer concordance invariants differ \cite{sneqtau}.
Or to mention an extreme example,
see Freedman et al. \cite{manandmachine} for a link with 222 crossings 
whose Rasmussen invariant was
worthwhile to compute at the time, because it could have provided a counterexample to
the four-dimensional smooth Poincar\'{e} conjecture (this approach was later rebuted
by Akbulut \cite{akbulut}).

Bar-Natan's extension of $\mathfrak{sl}_2$--homology from link diagrams to
tangles \cite{tanglescobordisms} (see also Khovanov \cite{khtangles}) led subsequently
to a divide-and-conquer algorithm to compute $\mathfrak{sl}_2$--homology \cite{fastcompu}.
The speed of this algorithm depends primarily on the \emph{girth} of the link diagram:
this is the maximal number of intersection points of a horizontal line with the diagram
(see e.g. \cite{freedman}, and cf. section \ref{sec:algo} for details).
An implementation by Green and Morrison called \texttt{JavaKh} \cite{javakh}
is able to compute the homology of knots of girth up to 14, e.g. the $(8,7)$-torus knot. 
%
Mackaay and Vaz \cite{vazuniverse} and Morrison and Nieh \cite{su} then extended
$\mathfrak{sl}_3$--homology to tangles, and the latter describe in detail the ensuing algorithm.
In this text, we present an implementation of this algorithm as a C++-program called \texttt{FoamHo} \cite{foamho}.

The algorithm can be improved by gluing sub-tangles in a more flexible way,
along a \emph{sub-tangle tree} 
instead of one after the other.
This leads to the notion of the \emph{recursive girth} of a link,
which replaces the girth as main factor limiting calculation speed.
Even though this improvement is not implemented in \texttt{FoamHo} yet, the program is still fast enough
to calculate the integral homology, reduced or unreduced, of links with girth up to 10, such as the $(6,5)$--torus knot.
The $\mathfrak{sl}_3$--concordance invariant $s_3$, as defined by Lobb
\cite{lobbconcordance} (see also Wu \cite{wu3} and Lobb \cite{lobbonroz}),
may (in most cases) be extracted from the $\mathfrak{sl}_3$--homology
by means of the spectral sequence converging to the filtered version of homology.
This method was used for $\mathfrak{sl}_2$--homology by Freedman et al. \cite{manandmachine}.
It does not depend on the conjectured convergence of the spectral sequence on the second page.

The most striking calculatory result obtained with \foamho{}
concerns this $s_3$--invariant, which is the $\mathfrak{sl}_3$--analogue of the Rasmussen invariant $s_2$.
Those two invariants share many properties, e.g. they agree on homogeneous and quasi-positive knots,
and until now it was not known whether they were actually equal or not (see Lobb \cite[Conjectures 1.5, 1.6]{lobbconcordance}).
\begin{figure}
\myhfill
\afig[height=3.81cm]{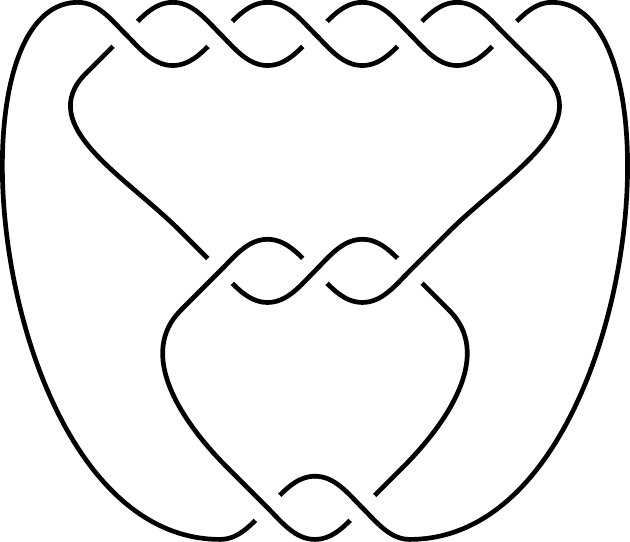}
\myhfill
\includegraphics[height=3.81cm]{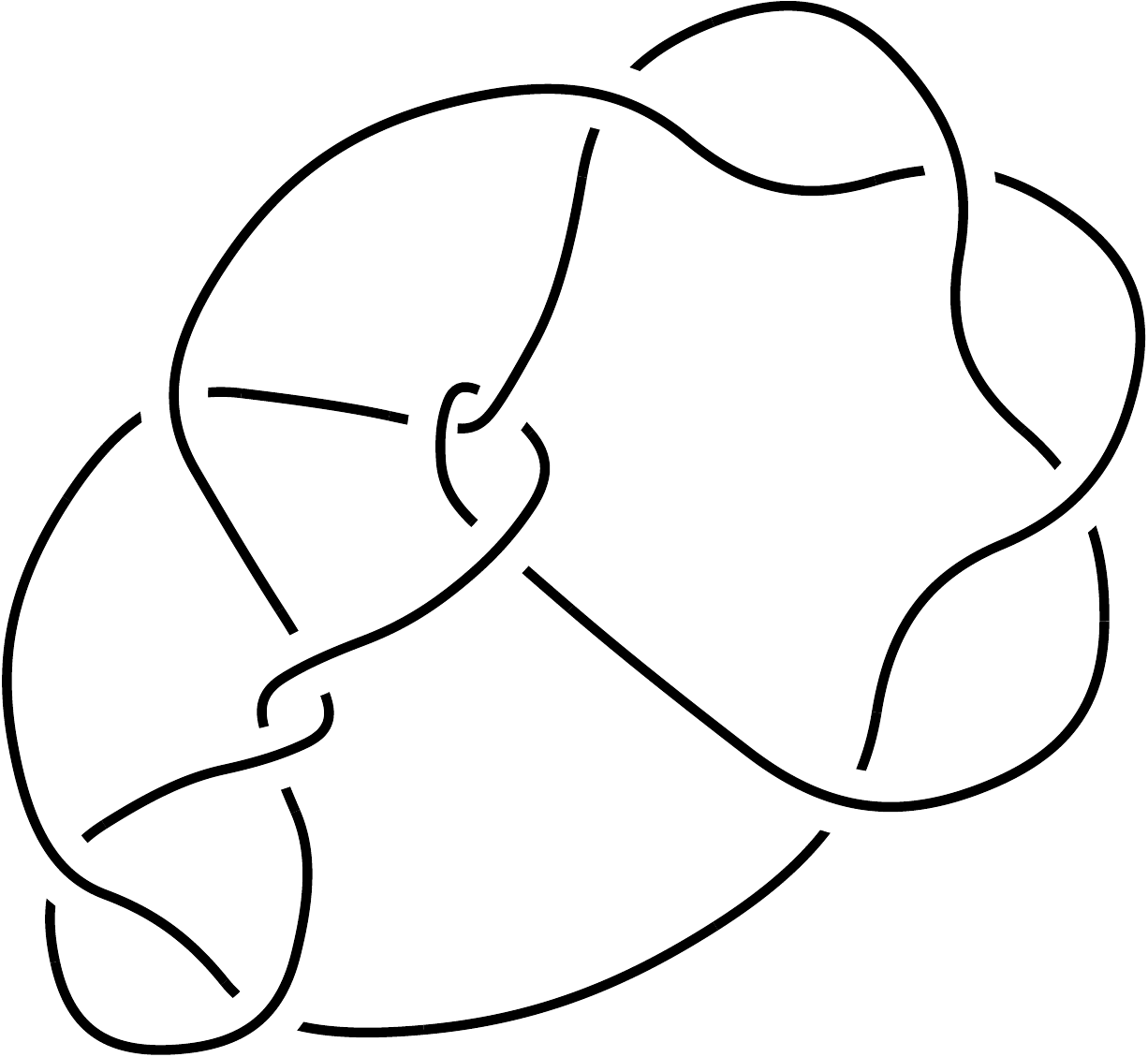}
\myhfill
\caption{On the left, the $(5,-3,2)$--pretzel knot (aka $10_{125}$), the smallest knot with $s_3 \neq s_2$;
on the right (drawn with knotscape \cite{knotscape}), 12n$_{340}$, the smallest knot with $|s_3| > |s_2|$.}
\end{figure}
\begin{conj} Let $s_3$ be the $\mathfrak{sl}_3$--concordance invariant, normalised to take 
the same value as the Rasmussen invariant on the trefoil.
If $a > b \geq 3$, $c \geq 2$, and $a + 1 \equiv b + 1 \equiv c \equiv 0 \pmod{2}$,
then the $(a,-b,c)$--pretzel knot $P(a,-b,c)$ has $s_3$--invariant
\[
s_3(P(a, -b, c)) = a - b + \delta_c,\]
where $\delta_2 = -1$ and $\delta_c =-2$ for $c > 2$.
\end{conj}
This statement is called a ``conjecture'' since it is established by \foamho-calculations
only for small values of $a$, $b$ and $c$. In addition, however, we have a prove for the case $c > 2$
which does not rely on computer calculations \cite{lewarkprepare}.
See section \ref{mainproof} for further details.

If one appends $b > c$ to the hypotheses of the conjecture,
then the $(a,-b, c)$--pretzel knot is quasi-alternating (see Champanerkar and Kofman \cite{quasialtpretzels} and Greene \cite{greene}),
and hence its $s_2$--invariant equals its classical knot signature (see Manolescu and Ozsváth \cite{quasialt}):
$s_2(P(a, -b, c)) = \sigma(P(a, -b, c)) = a - b$.
This value of the signature can be easily computed using G\"oritz matrices and
the formula of Gordon and Litherland \cite{gordon}.
So for this infinite family of pretzel knots, the $s_2$-- and $s_3$--invariant differ,
and the latter gives a weaker bound for
the slice genus than the former. However, there are knots
for which this is different, e.g. $s_3(12n_{340}) = 1$ and $s_2(12n_{340}) = 0$.

The rest of the paper is organised as follows:
section \ref{sec:2} provides a brief, but essentially self-contained definition
of the $\mathfrak{sl}_3$--foam homology as defined by Khovanov \cite{khovanovsl3},
Mackaay and Vaz \cite{vazuniverse}, and Morrison and Nieh \cite{su}.
We also define reduced homology in the framework of foams.
In section \ref{sec:3.1}, we give an account of the divide-and-conquer algorithm,
including the computation of integral homology and the
acceleration of the algorithm using a sub-tangle tree.
In section \ref{sec:3.2}, we discuss how to extract $s_3$ from the $\mathfrak{sl}_3$-homology.
Section \ref{sec:3.3} addresses some particular implementation
issues and their resolution in \texttt{FoamHo}, and section \ref{sec:3.4} presents the usage and characteristics of the program itself.
Finally, section \ref{sec:3.5} states some results of \texttt{FoamHo} calculations,
and compares them against previously known results.

\textbf{\emph{Acknowledgements}}\quad
I would like to thank Pedro Vaz and Alexander Shumakovitch for encouraging me to
pursue this subject and write a program, Louis-Hadrien Robert for all the inspiring discussions,
and Christian Blanchet for his continuous support and advice.
Thanks to Nils Carqueville, Pedro Vaz, Scott Morrison and Christian Blanchet for
comments on the first version of the paper.

\section{The $\mathfrak{sl}_3$-foam homology of tangle diagrams}
\label{sec:2}
This section gives a definition of the $\mathfrak{sl}_3$--polynomial
and its categorification, the $\mathfrak{sl}_3$--foam homology.
Except for the use of a generalised definition of planar algebras (sec. \ref{sec:pa}) to formalise localness,
and the definition of reduced homology using foams (sec. \ref{sec:red}), this section
contains nothing essentially new.
We just review the parts of \cite{khovanovsl3,su} which are relevant
to the purpose of this section -- which is to provide
a self-contained definition of $\mathfrak{sl}_3$--homology
with the objective of calculation in mind.
Instead of choking tori we use dots, like Khovanov \cite{khovanovsl3}
and Mackaay and Vaz \cite{vazuniverse}.
The origin of webs and the $\mathfrak{sl}_N$--polynomials lie in
representation theory \cite{RT,kuperberg},
an aspect we we will not dwell on.
\subsection{The $\mathfrak{sl}_3$--polynomial, naively}
The $\mathfrak{sl}_3$--polynomial can be defined by a single skein relation
involving only link diagrams.
We will instead use the two skein relations (Sk$^+$) and (Sk$^-$), see below,
because this allows a categorification using foams.
These skein relations involve webs:
a \emph{closed web} is a plane oriented trivalent graph
whose every vertex is either a source or a sink,
and that may have vertex-less circles as additional edges.

A \emph{tangle diagram} is the generic intersection of a link diagram
with a disc; generic means that the disc's border intersects the diagram's
strands transversely, and does not pass through a crossing.
Let us define a map $V$ from the set of smooth isotopy classes of tangle diagrams
to the free \qRing--module on the set of smooth isotopy classes of closed webs.
The map $V$ is uniquely determined by the following two local relations,
which are interpreted naively for now (i.e. apply these relations to all
crossings of the link at once, then expand):
\begin{align}
\tag{Sk$^-$} V\left(\avcfig{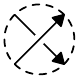}\right) & = q^2   \cdot V\left(\avcfig{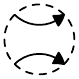}\right) - q^3   \cdot V\left(\avcfig{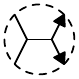}\right)\qquad\text{and} \\*
\tag{Sk$^+$} V\left(\avcfig{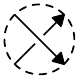}\right) & = q^{-2}\cdot V\left(\avcfig{03-Eqsign}\right) - q^{-3}\cdot V\left(\avcfig{02-H}\right).
\end{align}

Next, we define an evaluation \leftright{\langle}{\rangle} of closed webs,
called the \emph{Kuperberg bracket} \cite{kuperberg,jaeger},
which associates to a closed web a Laurent polynomial in $q$.
This evaluation is given by the four relations
\savebox\strutbox{$\vphantom{ V\left(\avcfig{32-neg}\right) = q^2   \cdot V\left(\avcfig{03-Eqsign}\right) - q^3   \cdot V\left(\avcfig{02-H}\right).}$}
\begin{align}
\tag{C} \left\langle\avcfig{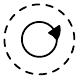}\right\rangle = 
\left\langle\avcfig{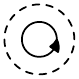}\right\rangle & = (q^{-2} + 1 + q^2)\cdot\left\langle\avcfig{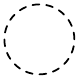}\right\rangle,\\*
\tag{D} \left\langle\avcfig{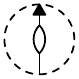}\right\rangle & = (q^{-1} + q)\cdot \left\langle \avcfig{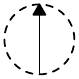}\right\rangle,\\*
\tag{S} \left\langle\avcfig{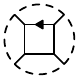}\right\rangle & = \left\langle \avcfig{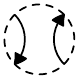}\right\rangle + \left\langle \avcfig{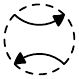}\right\rangle\qquad\text{and}\\*
\tag{U} \langle W_1 \sqcup W_2\rangle & = \langle W_1\rangle\cdot \langle W_2\rangle.
\end{align}
The $\mathfrak{sl}_3$--polynomial, which associates a Laurent polynomial in $q$
to a link diagram, can now be obtained by composing $V$ with \leftright{\langle}{\rangle},
and identifying the empty web with $1$.
Categorifying this construction is going to yield the $\mathfrak{sl}_3$--homology.
However, it is advantageous to formalise the \emph{localness} of the relations (Sk$^{\pm}$, C, D, S) before proceeding.
\subsection{Planar algebras}\label{sec:pa}
While 2--categories do give a framework for webs and foams, they make sense only if one aims at
interpreting webs as maps between oriented 0--manifolds; this aspect is not essential
to the calculation of $\mathfrak{sl}_3$--homology, and so we use planar algebras instead.
Planar algebras were introduced by Jones \cite{jones} to identify subfactors.
They were subsequently used to describe locally defined knot invariants such as
the Jones polynomial. Bar-Natan \cite{tanglescobordisms} introduced
a categorified version of planar algebras called \emph{canopolis} to describe
Khovanov homology, a method adaptable to $\mathfrak{sl}_3$--homology \cite{su}.
We will use a slightly generalised version of planar algebras,
working over arbitrary monoidal categories
instead of over the category of vector spaces over a fixed field.
In this way, a canopolis is a planar algebra as well.

Let $B_0 \subset \mathbb{R}^2$ be a closed disc, and $B_1, \ldots B_n \subset B_0^{\circ}$ be $n$ pairwisely
disjoint smaller closed discs.
Punching out these discs yields a disc with holes,
$\holy{H} = B_0 \setminus \bigcup_{i=1}^n B_i^\circ$.
Let $M$ be a compact oriented one-dimensional smooth submanifold of $\holy{H}$ with $M \cap \partial \holy{H} = \partial M$;
in other words, $M$ is a collection of circles and of intervals whose endpoints lie on the boundary of the big discs
or one of the smaller discs.  An \emph{input diagram} consists of $M$, $\holy{H}$, and on each boundary component of $\holy{H}$
a base point which is not in $\partial M$. We consider input diagrams up to smooth
isotopy, in the course of which boundary points of $M$ may not cross the base points.

For every $i \in\{0,\ldots n\}$, the intersection $M\cap B_i$ is a finite set;
at each of its points, the corresponding interval of $M$ is either oriented
towards the boundary ($+$ for $i > 0$, $-$ for $i = 0$), or away from it ($-$ for $i > 0$, $+$ for $i = 0$).
Moreover, these points have a canonical order, given by starting from the base point
and walking once around the circle in the counterclockwise direction.
Thus the isotopy type of
$M\cap B_i$ may be
written as a sign-word $\varepsilon_i$, i.e. a word over the alphabet $\{+,-\}$.
The \emph{boundary} of $\holy{H}$ is the tuple $(\varepsilon_0, \ldots \varepsilon_n)$.
\begin{figure}
\hfill\afig{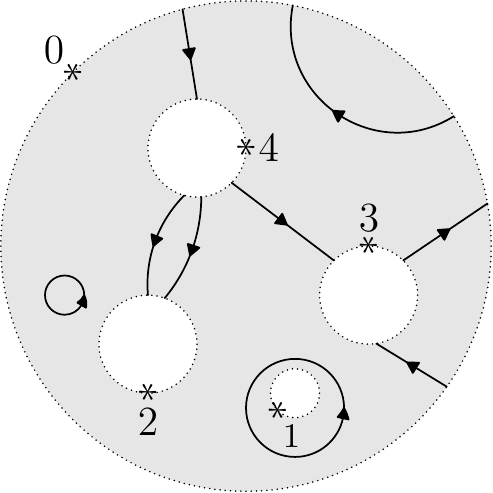}\hfill\rule{0pt}{0pt}
\caption{Example of an input diagram (also called spaghetti-and-meatballs diagram)
with boundary (+\,--\,+\,--\,+, $\varnothing$, ++, ++\,--\,, +\,--\,\,--\,\,--\,).
}
\label{fig:spaghetti}
\end{figure}

Now suppose $(M, \holy{H})$ and $(M', \holy{H}')$ are two input diagrams, such that $\varepsilon_0' = \varepsilon_i$
for a fixed $k \in \{1,\ldots n\}$. Then $(M', \holy{H}')$ may be
shrunk and glued into $B_k$, base point on base point and boundary points on boundary points,
resulting in a new input diagram with $n + n' - 1$ holes. 

Let $I$ be a subset of the set of all sign-words.
Let $\category{C}$ be a monoidal category, in the easiest case just the category $\category{Set}$ of sets,
and in the classical case the category of vector spaces over a fixed field.
Then a \emph{planar algebra} $\pa{P}$ over $I$ and $\category{C}$ consists of the following data:
\begin{itemize}
\item For each $\varepsilon \in I$, an object $\pa{P}_{\varepsilon} \in \category{C}$.
\item For each input diagram $\holy{H}$ with boundary $(\varepsilon_0, \ldots \varepsilon_n)$ such that $\forall i: \varepsilon_i \in I$,
a \category{C}--morphism
\[
\pa{P}_{\holy{H}}: \bigotimes_{i=1}^n \pa{P}_{\varepsilon_i} \to \pa{P}_{\varepsilon_0}.
\]
\end{itemize}
This data is required to satisfy the following axioms:
\begin{itemize}
\item Suppose $\holy{H}$ is an input diagram with $n = 1$ and $\varepsilon_0 = \varepsilon_1$
that consists only of appropriately oriented radial strands.
Then $\pa{P}_\holy{H}: \pa{P}_{\varepsilon_0} \to \pa{P}_{\varepsilon_1}$ is the identity morphism.
\item Let $\holy{H}$ and $\holy{H}'$ be two input diagrams with boundary $(\varepsilon_0, \ldots \varepsilon_n)$ and
$(\varepsilon'_0, \ldots \varepsilon'_n)$, respectively.
Suppose that for a fixed $k\in \{1,\ldots n\}$, $\varepsilon'_0 = \varepsilon_k$. Let $\holy{H}''$ be
the input diagram obtained from gluing $\holy{H}'$ into the $k$--th hole of $\holy{H}$.
Then the morphism $\pa{P}_{\holy{H}''}$ is equal to the composition of the morphisms $\pa{P}_{\holy{H}}$ and $\pa{P}_{\holy{H}'}$, i.e.
\[
\pa{P}_{\holy{H}''} = \pa{P}_\holy{H} \circ \bigl(\id_{\bigotimes_{i=1}^{k-1} \pa{P}_{\varepsilon_i}} \otimes \pa{P}_{\holy{H}'} \otimes 
\id_{\bigotimes_{i=k + 1}^{n} \pa{P}_{\varepsilon_i}}\bigr).
\]
\end{itemize}
\ \\

If $F: \category{C} \to \category{C}'$ is a monoidal functor, one may define the planar
algebra $F(\pa{P})$ over $I$ and $\category{C}'$ by $F(\pa{P})_{\varepsilon} = F(\pa{P}_{\varepsilon})$
and
\[
F(\pa{P})_{\holy{H}}: \bigotimes_{i=1}^n F(\pa{P})_{\varepsilon_i} \to F(\pa{P})_{\varepsilon_0}
\]
to be the composition $F(\pa{P}_{\holy{H}}) \circ \gamma$, where $\gamma$ is the natural transformation
\[
\bigotimes_{i=1}^n F(\pa{P}_{\varepsilon_i}) \to F\Bigl(\bigotimes_{i=1}^n \pa{P}_{\varepsilon_i}\Bigr)
\]
which comes with the functor $F$ because it is monoidal.
Examples of this construction include, for a planar algebra $\pa{P}$ over $\category{Set}$,
replacing for all $\varepsilon \in I$ the set $\pa{P}_{\varepsilon}$
by the free $R$--modules for some ring $R$ by means of applying the left-adjoint of the
forgetful functor from the category of $R$--modules to $\category{Set}$;
or the quotient $\pa{P}$ by an equivalence relation on $\pa{P}$, by which we mean a collection
of equivalence relations on all the $\pa{P}_{\varepsilon}$ which respect the planar
algebra structure.

Suppose $\pa{P}$ and $\pa{P}'$ are planar algebras over $I$, $I'$ and $\category{C}$, $\category{C}'$,
respectively, such that $I \subset I'$. A \emph{planar algebra morphism} from $\pa{P}$ to
$\pa{P}'$ consists of a functor $F: \category{C}' \to \category{C}$ and an $I$--indexed collection of
$\category{C}$--morphisms $\pa{P}_{\varepsilon} \to F(\pa{P}'_{\varepsilon})$
which respect the planar algebra structure, i.e. commute with the maps $\pa{P}_{\holy{H}}$ and
$\pa{P}'_{\holy{H}}$. The functor $F$ will typically be a forgetful functor.\\

Tangle diagrams with a base point on the boundary, considered -- as input diagrams -- up to smooth isotopy,
form a planar algebra $\pa{T}$ over $\category{Set}$ and the
set $I_0$ of sign-words $\varepsilon_1 \cdots \varepsilon_m$ with $\sum_{j=1}^m \varepsilon_j = 0$.
Let us elaborate this example:
the planar algebra $\pa{T}$ associates to a sign-word $\varepsilon$ with an equal number of both signs
the (countably infinite) set $\pa{T}_{\varepsilon}$ of all tangles diagrams with boundary $\varepsilon$, modulo
smooth isotopy; and to
an input diagram $H$ with $n$ holes (see e.g. fig. \ref{fig:spaghetti}) a function which maps a tuple
of $n$ tangle diagrams with appropriate boundaries to a new, bigger tangle diagram, by gluing each of the $n$
tangle diagrams into the corresponding hole of $H$. One easily verifies that the planar algebra axioms are satisfied.
\subsection{The $\mathfrak{sl}_3$--polynomial in the context of planar algebras}
Suppose the unit circle intersects a closed web generically; as in the definition of
tangle diagrams, this means that the circle intersects the edges of the closed web transversely,
and does not pass through a vertex.
Then the intersection of the closed web
with the unit disc is called a \emph{web}.
As for input diagrams, we fix a base point on the boundary of a web,
and encode the isotopy type of the boundary by a sign-word, in which 
$+$ stands for a strand oriented away from the boundary, $-$ for a strand oriented towards it.
Note that a sign-word $\varepsilon_1\cdots \varepsilon_m$ is the boundary of some web
if and only if $\sum_{j=1}^m \varepsilon_j \equiv 0 \pmod{3}$.
Denote by $I_3$ the set of all such sign-words.
Webs, up to smooth isotopy,
form a planar algebra $\pa{W}$ over $I_3$ and $\category{Set}$.

Let $\pa{W}^q_{\varepsilon}$ be the free $\qRing$--module on $\pa{W}_{\varepsilon}$.
Then $\pa{W}^q$ forms a planar algebra over $I_3$ and the category of $\qRing$--modules.
In $\pa{W}^q$, we may interpret
the relations (C), (D) and (S) as relations on $\pa{W}^q_{\varnothing}$,
$\pa{W}^q_{-+}$ or $\pa{W}^q_{+-}$ and $\pa{W}^q_{-+-+}$ or $\pa{W}^q_{+-+-}$, respectively.
Denote by $\pa{W}^{qr}$ the quotient by the generated equivalence relation.
In this context, the relation (U) is implied by the compatibility of the equivalence relation with
the planar algebra structure.

\label{rmk:webdiffeo}
Let $W, W'$ be two webs and $\varphi: W\to W'$ a diffeomorphism -- just of the
webs themselves, not taking into account the ambient discs.
We call $\varphi$ a \emph{web diffeomorphism} if it preserves the order of
the boundary points, and the cyclic ordering of edges around vertices.
Note that in the quotient $\pa{W}^{qr}$, the equivalence class of a web is already determined
by its \emph{web diffeomorphism} type. In $\pa{W}$, this distinction is slightly coarser than the
isotopy type, since e.g. web diffeomorphisms do not take the
orientation and relative position of closed components into account.

The two skein relations (Sk$^{\pm}$) determine a unique
morphism $V: \pa{T} \to \pa{W}^q$ of planar algebras, $\pa{T}$
being the planar algebra of tangles.
A link diagram $L$ may be seen as element of $\pa{T}_{\varnothing}$.
The equivalence class $[V(L)] \in \pa{W}^{qr}$
has a unique member that is a $\qRing$--multiple of the empty web.
The coefficient equals the $\mathfrak{sl}_3$--polynomial of the link diagram.
Reidemeister invariance may be shown by proving that the tangle diagrams
with two, four and six boundary points corresponding to the Reidemeister
moves I, II and III have in each case the same image under $V$.
\subsection{The $\mathfrak{sl}_3$--homology in the context of a special kind of planar algebras: canopolis}
\label{sec:can}
To categorify the $\mathfrak{sl}_3$--polynomial, one needs to understand foams, the cobordisms of webs.
Suppose that for all $i \in \{1,2,3\}$, $\Sigma^i$ are compact oriented smooth
(generally not connected) 2--manifolds
with $m$ boundary components $S^i_1, \ldots, S^i_m$ each.
Let $\varphi_j: S^1_j \to S^2_j$ and $\psi_j: S^1_j \to S^3_j$
be orientation preserving diffeomorphisms.
Consider the quotient of $\Sigma^1 \sqcup \Sigma^2 \sqcup \Sigma^3$
by the equivalence relation generated by $[x] = [\varphi_j(x)] = [\psi_j(x)]$ for all $j$.
The images of the $S^1_j$ in the quotient are called \emph{singular circles},
and the images of connected components of the $\Sigma^i$
are called \emph{facets}. There are three facets adjacent to
each singular circle.
Associate a non-negative integer to each facet. Such an integer $d$
will graphically be represented by drawing $d$ dots on the facet,
which may roam the facet freely, but may not cross a singular circle.
Such a quotient, together with the dots and with a choice of cyclic ordering
of the three facets around each singular circle is called a \emph{prefoam}.

\begin{figure}
\hfill\afig{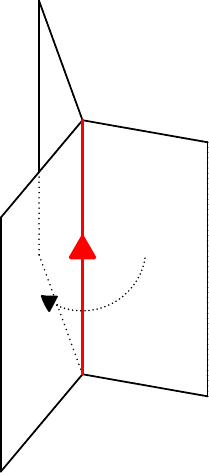}\hfill\rule{0pt}{0pt}
\caption{Cyclic ordering of facets around a singular circle of a closed foam.
Singular circles are drawn as a thick red line.}
\label{fig:cyclicord}
\end{figure}

Now consider a smooth embedding of a prefoam into $\mathbb{R}^3$,
i.e. an embedding that is smooth on the facets and on the singular circles.
Such an embedding induces cyclic orderings of the facets around each singular circle by
the left-hand rule (see figure \ref{fig:cyclicord}).
If these cyclic orderings agree with those given by the prefoam,
the image of the embedding is called a \emph{closed foam}.
Under the following conditions, 
the intersection of a closed foam with the cylinder $B \times [0,1] = \{(x,y,z) \mid x^2 + y^2 \leq 1\text{ and } 0 \leq z \leq 1\}$
is called a \emph{foam}:
\begin{itemize}
\item The boundary of the cylinder intersects the facets and singular circles of the closed foam transversely.
\item The side $(\partial B) \times [0,1]$ of the cylinder intersects the closed foam in finitely many vertical lines,
and is disjoint from all singular circles.
\item The intersections with the top and bottom of the cylinder are webs. 
\item The base point of the top and the base point of the bottom web have the same $x$-- and $y$--coordinates.
\end{itemize}
We consider foams up to isotopies which, on the side of the cylinder, do not
depend on the $z$--coordinate.
A connected component of the intersection of a singular circle with the cylinder is called a \emph{singular edge}.
The tangle diagrams on the bottom of the cylinder is called \emph{domain} of the foam,
and the \emph{codomain} of the foam is defined as the tangle diagram on the top
of the cylinder, with the orientation of each edge reversed.
As usual, let $\chi$ denote the Euler characteristic. Then the \emph{degree} of a foam $f$ is defined by
\[
\deg f = 
\chi(\text{domain of }f)
+ \chi(\text{codomain of }f)
+ 2(\text{total number of dots on }f)
- 2\chi(f).
\]

Foams can be glued in two ways:
if the domain of one foam agrees with the codomain of another, by stacking them on top of each other.
Or, by gluing them into the cylindrical holes of a thickened input diagram.
The degree is additive with respect to both of these operations.

Webs with a fixed boundary and the foams between them thus constitute a \emph{graded category},
i.e. a category whose morphisms have an integral rank
which is additive under composition.
Let us define a planar algebra
$\pa{W}^c$ over $I_3$ and the category $\category{GCat}$
of small graded categories: \footnote{The $c$ superscript stands for categorification.}
to $\varepsilon\in I_3$, associate the category whose set of objects is
$\pa{W}_{\varepsilon}$, and whose morphisms $W \to W'$ between two webs
$W, W' \in \pa{W}_{\varepsilon}$ are the foams with domain $W$ and codomain $W'$.
If $\holy{H}$ is a planar input diagram, then $\pa{W}^c_{\holy{H}}$ is the functor
that acts as $\pa{W}_{\holy{H}}$ on the objects,
and glues foams into a thickened version of $\holy{H}$.

Next, $\pa{W}^{cq}$ may be constructed by applying a functor from $\category{GCat}$
to $\category{ACat}$, the category of small additive categories:
replace webs by $\qRing$--linear
combinations of webs, and foams by matrices of $\mathbb{Z}$--linear combination of foams,
where morphisms from $q^{\alpha}\cdot W$ to $q^{\beta}\cdot W'$ are the foams with degree $\beta - \alpha$.
So the categories $\pa{W}^{cq}_{\varepsilon}$ are not graded, but have instead a shift operator
for their objects.
In this planar algebra, consider the following morphisms:
\begin{align}
\tag{C$^c$} 
\avcfig{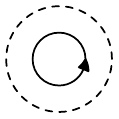} &
\xymatrix@=5cm@M=3mm{
\ar@<0.3em>[r]^-{\left(-\:\avcfig{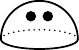}\quad -\:\avcfig{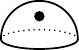}\quad -\:\avcfig{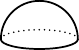} \right)^{\top}}   &
\ar@<0.3em>[l]^-{\left(\avcfig{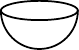}\quad \avcfig{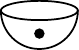}\quad \avcfig{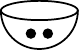} \right)}
} (q^{-2} + 1 + q^2)\cdot\avcfig{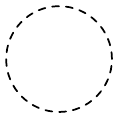},\\[5mm]
\tag{D$^c$} 
\avcfig{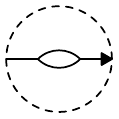} &
\xymatrix@=5cm@M=3mm{
\ar@<0.3em>[r]^-{\left(\:-\:\avcfig{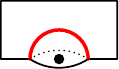}\quad \avcfig{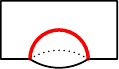}\: \right)^{\top}} &
\ar@<0.3em>[l]^-{\left(\:\avcfig{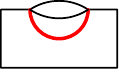}\quad \avcfig{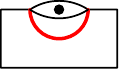}\: \right)}
} (q^{-1} + q)\cdot\avcfig{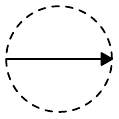}\qquad\text{and}\\[5mm]
\tag{S$^c$} 
\avcfig{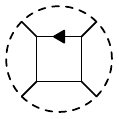} &
\xymatrix@=5cm@M=3mm{
\ar@<0.3em>[r]^-{\left(\:-\:\avcfig{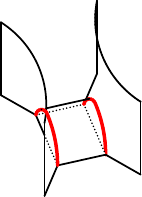}\quad -\:\avcfig{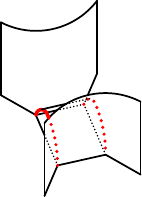}\: \right)^{\top}} &
\ar@<0.3em>[l]^-{\left(\:\avcfig{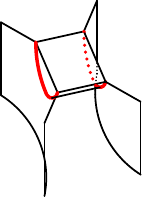}\quad \avcfig{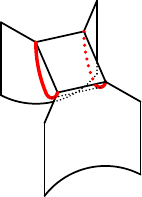}\: \right)}
} \avcfig{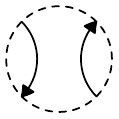} + \avcfig{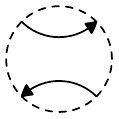}
\end{align}
Demanding that these three pairs of morphisms are mutually inverse
 -- for any placement of the base point --
generates an equivalence relation on the planar algebra $\pa{W}^{cq}$.
For example, composing the two morphisms in (C$^c$) one way yields the so-called
surgery relation, and composing them the other way yields the evaluation of the sphere
with four or less dots.

Let $\pa{W}^{cqr}$ be the quotient of $\pa{W}^{cq}$ by this equivalence relation.
It is a consequence of \cite{khovanovsl3} that the relation set
(C$^c$), (D$^c$), (S$^c$) is equivalent to the original relation set given
by Khovanov. Thus the elegant, but non-constructive universal BHMV-construction \cite{bhmv}
may be circumvented.

Let us consider some examples of $\pa{W}^{cqr}_{\varepsilon}$.
For $\varepsilon = \varnothing$, the isomorphism classes of objects of this category are in correspondence
with the elements of the free $\qRing$--module generated by the empty web $\varnothing$.
Non-zero morphisms $q^{\alpha}\cdot\varnothing \to q^{\beta}\cdot\varnothing$ exist only if $\alpha = \beta$, and in this
case the morphism $\mathbb{Z}$--module is just $\mathbb{Z}$.
The $\qRing$--module Khovanov \cite{khovanovsl3} associates to a closed web $W$ can be recovered
as
\[
\bigoplus_{\alpha\in\mathbb{Z}} \Hom_{\pa{W}^{cqr}_{\varnothing}}(q^0\cdot \varnothing, q^{\alpha} \cdot W).
\]

Similarly, the isomorphism classes of  objects of $\pa{W}^{cqr}_{+-}$ are in correspondence with the
elements of the free $\qRing$--module generated by the web \avcfig[scale=0.75]{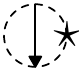}, which is just an interval.
The morphisms are more complicated, however, since there are non-zero foams of three
different degrees: rectangles with none, one, or two dots.
So the morphism $\mathbb{Z}$--module from $q^{\alpha}\cdot\avcfig[scale=0.75]{100-straightUpSmall} \longrightarrow q^{\beta}\cdot\avcfig[scale=0.75]{100-straightUpSmall}$
is $\mathbb{Z}$ if $\beta-\alpha \in \{0, 2, 4\}$, and trivial otherwise.

\label{rmk:foamdiffeo}
Let a diffeomorphism between two foams
be called a \emph{foam diffeomorphism} if its restriction to the top (bottom)
of the cylinder constitutes a web diffeomorphisms between the domains (codomains)
of the two foams. In $\pa{W}^{cqr}$, the equivalence class of a foam is already
determined by its foam diffeomorphism type; this is not the case in $\pa{W}^{cq}$,
where e.g. a cylindrical foam tied into a knot is not the identity of the circular web.

Finally, let $\pa{W}^{cqrt}$ be the planar algebra obtained from $\pa{W}^{cqr}$ by
setting $\pa{W}^{cqrt}_{\varepsilon}$ to be the category of bounded chain complexes (up to chain homotopy)
over $\pa{W}^{cqr}_{\varepsilon}$.
More formally, this means applying a monoidal functor $K: \category{ACat} \to \category{ACat}$.
The natural transformation
$K(C_1)\otimes K(C_2) \to K(C_1 \otimes C_2)$ is given by
\begin{equation}\tag{*}\label{tag:ccproduct}
(P_i, g_i)_i \otimes (Q_j, h_j)_j \mapsto
\Bigl(
\bigoplus_{i + j = k} P_i \otimes Q_j, \sum_{i + j = k} g_i \otimes \id_{Q_j} + (-1)^i \id_{P_i} \otimes h_j
\Bigr)_k.
\end{equation}
In the notation of complexes, the module at $t$--degree 0 is underlined.
We may now define the $\mathfrak{sl}_3$--chain complex as the planar algebra morphism
$V^c: \pa{T} \to \pa{W}^{cqrt}$ uniquely determined by the skein relations
\begin{align}
\tag{Sk$^{-c}$} V^c\left(\avcfig{32-neg}\right) & = \underline{q^2 \avcfig{03-Eqsign}} \xrightarrow{\avcfig{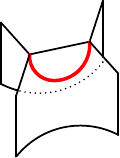}} q^3\avcfig{02-H}
\qquad\text{and}\\*
\tag{Sk$^{+c}$}
V^c\left( \avcfig{35-pos}\right) & = q^{-3} \avcfig{02-H} \xrightarrow{\avcfig{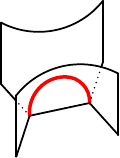}} \underline{q^{-2}\avcfig{03-Eqsign}}.
\end{align}
Note that these relations each just give the complex of one tangle diagram,
and no cone or flattening is involved. The minus-sign which must intervene
to arrange anti-commutativity of differentials along squares is hidden in 
(\ref{tag:ccproduct}). Although the placement of this sign is somewhat arbitrary,
$V^c$ is still unique, since chain complexes are considered only up to homotopy.

Reidemeister invariance may, just as for the $\mathfrak{sl}_3$--polynomial,
be shown by inspecting the small tangle diagrams corresponding to Reidemeister moves \cite[section 4.2]{su}.
So the $\mathfrak{sl}_3$--chain complex is, up to chain homotopy, an invariant of tangles.

A link diagram lives in $\pa{T}_{\varnothing}$, and is mapped to a chain complex in the
category of $\qRing$--linear combination of closed webs and matrices of $\mathbb{Z}$--linear
combination of foams between them. Due to the relations, there is a homotopy to a chain
complex which contains only empty webs and closed foams. The latter are in turn just integral
multiples of the empty foam. Setting the empty foam to $1$, a chain complex in the
category of free graded abelian groups is obtained. Its homology is the \emph{$\mathfrak{sl}_3$--homology}
of the link.
\subsection{Reduced $\mathfrak{sl}_3$--homology}\label{sec:red}
A \emph{reduced} version of $\mathfrak{sl}_N$--homology
has been introduced by Khovanov and Rozansky \cite{roz}
in the context of matrix factorisations.
This section contains a definition in the context of foams.

Let $D$ be a diagram of a link $L$ with a marked component. Cutting $D$ open at an arbitrary point on the marked component,
one obtains a tangle diagram $D'$ with boundary $+-$. It is well known that two
such tangle diagrams represent the same link with a marked component if and only if they are connected
by a finite sequence of Reidemeister moves. So, the homotopy type of the $\mathfrak{sl}_3$--chain
complex $C$ of $D'$ is an invariant of links with a marked component; in particular, it
is a knot invariant.
%
As noted in section \ref{sec:can}, if $C$ is fully simplified,
it contains only one kind of web -- an interval --, and three kinds
of foams: the identity foam of the interval with none, one or two dots.
Now map the identity foam to one,
and foams of higher degree to zero. Thus one obtains a chain complex
of free graded abelian groups. Its homology is the
reduced $\mathfrak{sl}_3$--homology of $L$,
an invariant of links with a marked component.
As simple examples show, its value may indeed depend on the choice of
marked component, and it does neither determine the unreduced version, nor
is it determined by it (see section \ref{sec:session}).

\section{Calculations}\label{sec:impl}
\subsection{The algorithm}\label{sec:algo}\label{sec:3.1}
The above definition of $\mathfrak{sl}_3$--homology gives a straightforward way of practical calculation:
for an $n$--crossing diagram, take $n$ copies of the chain complex of Sk$^{\pm c}$,
and take the tensor product (\ref{tag:ccproduct}) of these complexes as determined by the diagram.
This is the same as forming the \emph{cube of resolutions}.
Transform this chain complex of closed webs and foams between them
into a homotopic chain complex of empty webs and closed foams
using the three relations (C$^c$), (D$^c$), (S$^c$); finally, evaluate the closed foams,
and calculate the homology of the emerging integral chain complex.

In the algorithm described by Morrison and Nieh \cite{su}, it is just the order in which these steps are taken
which is changed: we do not apply (Sk$^{\pm c}$) to all crossings at once, but to one
after the other, and try to simplify the chain complex at each step as much as possible.
Examples for the manual application of this algorithm can be found in \cite{su}.

At each step, one manipulates a chain complex of tangles; if these tangles have less boundary points,
then there are fewer different tangles and cobordisms between them and
thus the calculations demand less memory and will, heuristically, go faster.
Therefore, one glues the crossings in such an order that the cardinal of the boundary
of the intermediate tangles is minimised (or at least, as low as one sees possible).
This minimum is precisely the \emph{girth} of the link diagram, which is thus the
main factor limiting the speed of the algorithm.

In this text, we propose a variation of this algorithm based on sub-tangle trees, which is potentially faster,
because intermediate tangles will have smaller boundary.
\begin{figure}
\hfill\afig{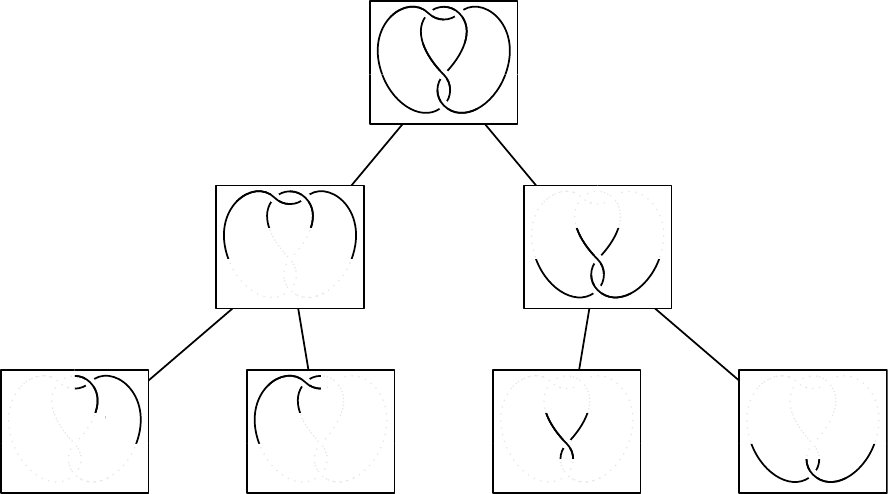}\hfill\rule{0pt}{0pt}
\caption{Example of a sub-tangle tree of the figure-eight-knot with girth 4.}
\label{fig:tree}
\end{figure}
\begin{defn}Let $D$ be a non-split link diagram with crossings enumerated from $1$ to $n$.
Decompose $D$ into small tangle diagrams $D_i$ that contain the $i$--th crossing, respectively.
A \emph{sub-tangle tree} of $D$ is a full binary tree with a tangle diagram at
each node, such that
\begin{itemize}
\item The leaves are decorated by the $D_i$.
\item The root is decorated by $D$.
\item Every node which is not a leaf has two children decorated with adjacent
tangle diagrams, and is decorated itself with the union of those two tangle diagrams.
\end{itemize}
Let the girth of a sub-tangle tree be the maximum number
of boundary points of the tangle diagrams at all nodes.
The \emph{recursive girth} of $D$ is the minimum of the girths of all its sub-tangle trees.
\end{defn}
In general, the recursive girth of a diagram is smaller than its girth. For example, pretzel
links have recursive girth 4, and girth 6; and the 222--crossing link considered in
\cite{manandmachine} has girth $\approx 24$, and recursive girth at most 16.

Let us give a description of the algorithm.
Given a link diagram $D$, find a sub-tangle tree of $D$ with small girth.
The algorithm consists in simplifying the $\mathfrak{sl}_3$--chain complex
of the tangle diagram at each node, one after the other. In the final step, the chain complex of the diagram
at the root is handled, which is $D$ itself. After simplifications, 
this will be a chain complex of vectors of $q$--shifted empty webs and matrices whose entries
are integer multiples of the empty foam; identifying the empty foam with $1$,
this becomes a complex of free graded abelian groups, and its homology may be calculated --
separately for each $q$--degree -- using the Smith normal form.

In the beginning, only the complexes at the leaves are known, which are given by
the relations (Sk$^{\pm c}$).
During each step, fix a node whose two children have complexes $C_1$ and $C_2$
that are already known and simplified.
Now, decorate this node with the tensor product
$C_1 \otimes C_2$ (see (\ref{tag:ccproduct})) and simplify this chain complex as described below.
After as many steps as $D$ has crossings, the process reaches the root and is finished.

There are two ways to simplify a chain complex:
firstly, apply the circle, digon or square relation wherever possible.
Secondly, apply Gau\ss' homotopy to all matrix entries which are
are plus or minus an invertible foam:
\begin{lem}[Gau\ss' homotopy, \cite{fastcompu}]
Over an additive category with an isomorphism $h$, the chain complex
\[
\xymatrix@=1.2cm{
\ldots \ar[r] & P \ar[r]^{\mymatrixtwov{*}{g}} &  Q\oplus R
\ar[r]^{\mymatrixfour{h}{i}{j}{k}} & S\oplus T
\ar[r]^{\mymatrixtwoh{*}{\ell}} & U \ar[r] & \ldots
}
\]
is homotopic to
\[
\xymatrix@=1.7cm{
\ldots \ar[r] & P \ar[r]^{\mymatrixsimple{g}} &  R
\ar[r]^{\mymatrixsimple{k - j \circ h^{-1} \circ i}} & T
\ar[r]^{\mymatrixsimple{\ell}} & U\ar[r] & \ldots.
}
\]
\end{lem}
The additive category in question is $\pa{W}^{cqr}_{\varepsilon}$,
whose objects are $R[q^{\pm 1}]$--linear combinations of webs with boundary $\varepsilon$,
and whose morphisms are matrices of $R$--linear combinations of foams.
The base ring $R$ is usually $\mathbb{Z}$;
but one may also just calculate homology over some field.
Gau\ss' homotopy may then be applied more often,
namely to all non-zero multiples of invertible foams instead of just plus or minus
an invertible foam. This may increase the algorithm's speed.

As a rule, use the circle/digon/square-relations only when it is not possible to
apply Gau\ss' homotopy. This is to prevent those relations from altering a foam
which is plus or minus the identity and could thus be removed.
\subsection{Extracting the $s_3$--concordance invariant from homology}\label{sec:rasextract}
\label{sec:3.2}
There is a spectral sequence $E_{\bullet}$ from the graded $\mathfrak{sl}_3$--homology
converging to a the filtered version.
Given the $\mathfrak{sl}_3$--homology of a knot, this allows in many cases to
determine its $s_3$--invariant.
The same method was used by Freedman et al.
for $\mathfrak{sl}_2$--homology \cite[section 5.2]{manandmachine}.

The filtered $\mathfrak{sl}_3$--homology has Poincar\'{e} polynomial
$q^{-2s_3}(q^{-2} + 1 + q^2)$. The differential on the $k$--th page of the spectral
sequence has $(t,q)$--degree $(1, -3k)$.
Let $KR_3$ be the Poincar\'{e} polynomial of the graded $\mathfrak{sl}_3$--homology of
a fixed knot. Then $E_{\bullet}$ gives rise to a decomposition of the form
\begin{equation*}
KR_3(t,q) = q^{-2s_3}(q^{-2} + 1 + q^2) + \sum_{k = 1}^\infty \zeta_k(t,q)\cdot (1 + tq^{3k}),
\end{equation*}
for some Laurent polynomials $\zeta_k(t,q)$ with non-negative coefficients.
Reversely, given the value of $KR_3$, one may easily determine all possible such decompositions.
The conjecture that $E_{\bullet}$ converges on
the second page translates as $\forall k \geq 2: \zeta_k(t,q) = 0$. Assuming the conjecture
is true, the above decomposition is unique; otherwise, however, it is generally not.
But even if there are are several possible decompositions, it may happen they all share
the same value for $s_3$. This need not be the case, either, as the example
$KR_3(t,q) = 1 + q^2 + q^4 + q^6 + tq^{12}$ demonstrates: either
$s_3 = -1, \zeta_1(t,q) = q^6, \zeta_{\neq 1}(t,q) = 0$, or $s_3 = -2, \zeta_2(t,q) = 1, \zeta_{\neq 2}(t,q) = 0$.
The first knots for which $s_3$ cannot be uniquely determined from $\mathfrak{sl}_3$--homology
are 12n$_{118}$, 12n$_{210}$, 12n$_{214}$ and 12n$_{318}$ (see section \ref{sec:session} for further examples).
At any rate, one obtains a list of possible values for $s_3$,
and for most small knots, only one value is possible.
\subsection{Implementation issues}
\label{sec:3.3}
\begin{figure}[h!]
\myhfill
\avcfig{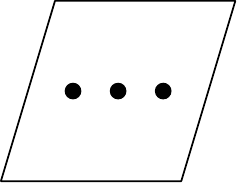} =\ 0
\myhfill
\avcfig{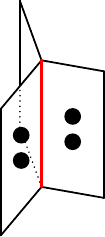}\quad =\quad 0
\myhfill
\avcfig{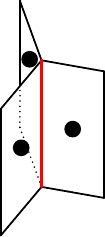}\quad =\quad 0
\myhfill
\\[4mm]
\myhfill
\avcfig{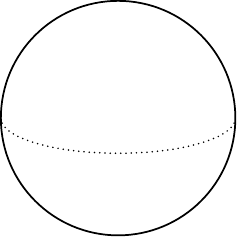}\quad =\quad \avcfig{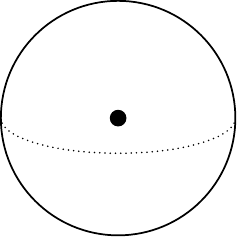}\quad =\quad 0
\myhfill
\avcfig{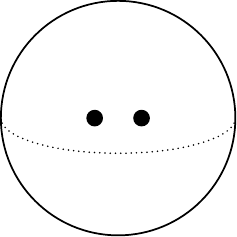}\quad =\quad $-1$
\myhfill\\[6mm]
\myhfill
\avcfig{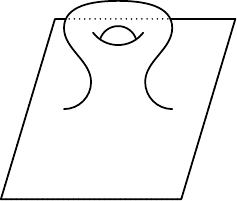} = $-3$ \avcfig{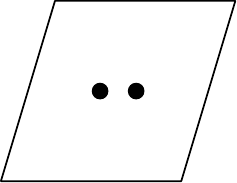}
\myhfill\\[4mm]
\begin{align*}
\avcfig{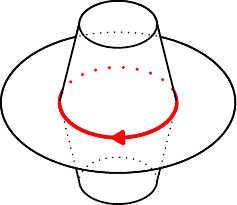}\quad =\quad
& \avcfig{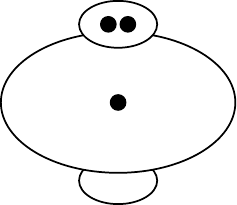}\ +\ 
\avcfig{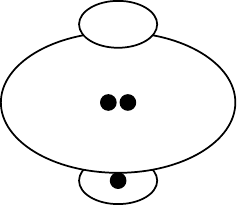}\ +\ 
\avcfig{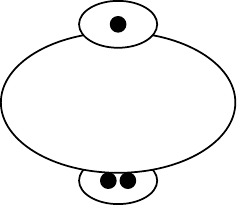} \\[2mm]
-\ 
& \avcfig{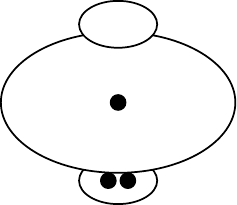}\ -\ 
\avcfig{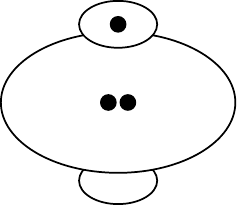}\ -\ 
\avcfig{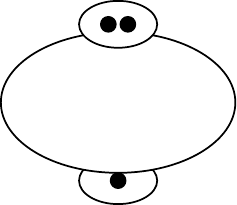}
\end{align*}
\caption{Foam relations}
\label{fig:foamrel}
\end{figure}

An implementation of the algorithm of section \ref{sec:algo} is possible because webs and foams
need only be considered up to web and foam diffeomorphisms (see sections
\ref{rmk:webdiffeo} and \ref{rmk:foamdiffeo}), and their
diffeomorphism type contains only finitely much information.
A web is determined by the following data:
\begin{itemize}
\item Its boundary, given as a sign-word.
\item The number of edges that are circles.
\item For each vertex, the triple of the vertices or boundary points
it is connected to, listed in counterclockwise order.
\end{itemize}
Likewise, a foam
is completely encoded by the following information:
\begin{itemize}
\item Its domain and codomain web.
\item The start and end point of every singular edge that is an interval.
\item For each singular edge, the three facets adjacent to it, in the order specified by the left-hand rule (see fig. \ref{fig:cyclicord}). 
\item The genus of and number of dots on each facet.
\item Each boundary component of each facet, given as an ordered tuple of edges;
in this tuple, edges of the domain and codomain web, and singular edges alternate.
\end{itemize}
In order to reduce the complexity of foams, \foamho{} applies the well-known
relations shown in fig. \ref{fig:foamrel} whenever possible. In particular,
these relations are sufficient to evaluate closed foams.
Gau\ss' homotopy is in fact only applied to plus or minus an identity foam;
there are other isomorphisms than those, but
they are difficult to detect for a computer program, and so
rare that looking for them does not seem worthwhile.
\subsection{\foamho, a $\mathfrak{sl}_3$--calculator}
\label{sec:3.4}
The algorithm described in the previous sections was implemented by the author
as a C++-program\footnote{Using the PARI/GP-library \cite{pari} to calculate the Smith normal form,
and the MPIR library \cite{mpir} for arbitrary precision integers and rationals.}.
It computes the integral or rational (homology over finite fields is not yet implemented),
reduced or unreduced homology of knot or link diagrams given in braid, planar diagram or Dowker-Thistlethwaite notation;
it also attempts to find the value of
the $s_3$--concordance invariant, using the spectral sequence from $\mathfrak{sl}_3$--homology
to filtered $\mathfrak{sl}_3$--homology.
If it is impossible to extract the $s_3$-invariant, a list of possible values is printed,
and the value which corresponds to the convergence of the spectral sequence on the second page is highlighted.

The current version does not make use of the sub-tangle trees yet,
and glues instead one crossing after the other.
Apart from that, there is surely some room for rendering the program faster and
less memory hungry by optimising the code, without changing the algorithm;
e.g., webs and foams are encoded in a redundant way, and tightening this would
reduce the memory consumption.

The program was baptised \foamho\ in recognition of Shumakovitch's \texttt{KhoHo} \cite{khoho}.
It has been released under the GPL\footnote{See \url{http://www.gnu.org/licenses/gpl.html}.},
and its source code as well as compiled versions for
Linux and Windows may be downloaded from \cite{foamho}.
There are also tables of the homology and the $s_3$--concordance invariant
of small knots and links available.

Let us give a rough idea of the capabilities of the program for large knots.
On an AMD Opteron (2.2 GHz), the (6,5)--torus knot's homology can be calculated
in six minutes, using 80 MB of RAM; the (8,5)--torus knot takes 50 minutes and
275 MB of RAM; but the (7,6)--torus knot is out of the reach of
5 GB of RAM. In other words, links of girth 10 can be calculated until ca. $40$
crossings, while links of girth 12 would demand a well-equipped computer.

Detailed usage instructions for \foamho\ 
can be found in the README-file which is distributed along with
the program's source code and binaries.
The following appetising example session demonstrates the computation of the
$(4,3)$--torus knot's homology.
User input begin with a dollar sign and is printed bold:\\
\begin{Verbatim}[formatcom=\bfseries\normalsize]
$ ./foamho -h
\end{Verbatim}
\begin{Verbatim}[formatcom=\normalsize]
foamho, a sl3-homology calculator, version 1.1.

Usage: foamho [OPTIONS] braid | pd | dt NOTATION

For example, the following three commands all compute the figure-8-knot's 
integral homology:

                   foamho braid aBaB
                   foamho pd "[[4,2,5,1],[8,6,1,5],[6,3,7,4],[2,7,3,8]]
                   foamho dt "[4,6,8,2]"

Options:
-q                 Compute rational homology instead of integral.
-r                 Compute reduced homology instead of unreduced. You may
                   give a number right after -r to indicate the marked
                   strand (useful for links).

-g                 Do not attempt to optimise the girth.
-v                 Display some progress information.
-vv                Display more detailed progress information.
-t                 Display time and memory consumption.
-h                 Display this help message and exit.

Written in 2012/2013 by Lukas Lewark, lewark@math.jussieu.fr.
All feedback is welcome.
\end{Verbatim}
\begin{Verbatim}[formatcom=\bfseries\normalsize]
$ ./foamho -t braid abababab    
\end{Verbatim}
\begin{Verbatim}[formatcom=\normalsize]
Girth-optimised link diagram (modified pd notation) [[2,4,3,1],[5,7,6,4],
    [6,9,8,3],[7,11,10,9],[10,13,12,8],[11,15,14,13],[14,16,1,12],
    [15,5,2,16]].
Girth: 6.
Calculating...
Done. Result:
Rational homology: (q^-14 + q^-12 + q^-10) + t^2(q^-16 + q^-14) +
    t^3(q^-22 + q^-20) + t^4(q^-20 + 2q^-18 + q^-16) +
    t^5(q^-26 + 2q^-24 + q^-22) + t^6(q^-22) + t^7(q^-28 + q^-26) +
    t^8(q^-32)
Total rank: 19
Rational homology is not self-dual => the link is chiral.
3-torsion: t^3(q^-18) + t^5(q^-22 + q^-20) + t^7(q^-26 + q^-24) +
    t^8(q^-30 + q^-28)
The 3-torsion part of homology is not self-dual => the link is chiral.
s_3-concordance invariant: -12
Run time in seconds: 2
Memory consumption in megabytes: 9.5
\end{Verbatim}
\subsection{Calculatory results}
\label{sec:3.5}
\label{mainproof}
\label{sec:session}
\myparagraph{Pretzel knots:}
Since pretzel knots have girth only 6, their homology can be computed particularly quickly.
In addition to \texttt{FoamHo} calculations, we use an inequality which
is stated by Lobb \cite{lobbineq} (cf. also Kawamura \cite{kawamura})
for the Rasmussen invariant, but holds for the $s_3$--invariant as well.
With these two tools, we can confirm the conjecture (see introduction) for the following values:
\begin{itemize}
\item $b \leq 49$, any $a > b$ and any $c \geq 4$ (an infinite family).
\item $b < a \leq 49$ and $c = 2$.
\end{itemize}
Furthermore, we are able to prove the conjecture for $c > 2$, and prove $\delta_2 \in \{-1,-2\}$.
The proof (paper in preparation \cite{lewarkprepare})
makes use of the tools developed by Rasmussen in \cite{somedifferentials}.
\myparagraph{$s_3$ and $s_2$ for small knots:}
Of the 59\,937 prime knots with 14 crossings or less, there are 361 for which
the $s_3$ and $s_2$ concordance invariants differ, and for 63
the absolute value of $s_3$ is greater, i.e. the lower bound to the slice
genus coming from $\mathfrak{sl}_3$ is better than the one coming from $\mathfrak{sl}_2$.
For all these 361 knots, the difference between $s_2$ and $s_3$ is $\pm 1$.
A 16--crossing prime knot example for a greater difference is given by
the conjecture (see introduction): the $(7,-5,4)$--pretzel knot,
with $s_3 = 0$ and $s_2 = 2$.

The $s_3$--concordance invariant cannot be determined for certain knots with
14 crossings or less, so the statistics in the preceding paragraph are based on the assumption that the
the spectral sequence to the filtered version of homology converges on the second page.
\myparagraph{$s_3$ and the Floer-concordance invariant $\tau$:}
Let $\tau$ be the knot concordance invariant coming from Floer homology, as defined by Ozsváth and Szabó \cite{os}.
Hedden and Ording \cite{sneqtau} found examples of knots $K$ for which $2 = s_2(K) \neq 2\tau(K) = 0$.
\foamho{} calculations yield the following results for these knots:
\begin{alignat*}{5}
s_3(D_+ (T_{2,3} , 2)) & = 2  &            s_3(D_+ (T_{2,7} , 8)) & \in \{2,3,4\} \\
s_3(D_+ (T_{2,5} , 5)) & \in \{2,3\}  &    s_3(D_+ (T_{2,7} , 7)) & \in \{2,3,4,5\} \\
s_3(D_+ (T_{2,5} , 4)) & \in \{2,3,4\}\qquad\qquad  &  s_3(D_+ (T_{2,7} , 6)) & \in \{2,3,4,5,6\}
\end{alignat*}
All but the first knot are examples of how the method of section \ref{sec:rasextract}
may fail to completely determine the value of $s_3$.
Under the assumption that the spectral sequence to the filtered version of homology
converges on the second page, however, each of the above knots $K$ satisfies $s_3(K) = 2$.
So these knots cannot be expected to give examples for $s_2 \neq s_3$,
but they demonstrate that $s_3 \neq 2\tau$.
\myparagraph{Thinness:}
Khovanov \cite{patterns} called a knot \emph{H-thin} if its unreduced Khovanov homology is supported in only
two diagonals, or, equivalently, if its reduced Khovanov homology is supported in only one.
Thinness was generalised to $\mathfrak{sl}_N$-homologies by Rasmussen \cite{rasbridge}.
Alternating knots are H-thin, but not necessarily \emph{$N$-thin} for $N>2$. 
%
%
We find the knots $11a_{263},12a_{36},12a_{694},12a_{804},12a_{811},12a_{817},12a_{829},12a_{832}$
to be the smallest examples of alternating knots that are not 3--thin.
None of these knots is two-bridge, which agrees with
Rasmussen's theorem \cite{somedifferentials} that two-bridge knots are $N$--thin for all $N$.
\myparagraph{Mutation:}
Integral reduced $\mathfrak{sl}_3$--homology has been proven to be invariant
under mutation of knots by Jaeger \cite{khrozmutation}; for unreduced homology, invariance appears to
be an open question. Our calculations confirm that
all mutant families up to 13 crossings have the same reduced and unreduced integral $\mathfrak{sl}_3$--homology.
We use the lists provided by Stoimenow \cite{stoimenow}.
\myparagraph{Torsion:}
All prime knots and links for which the homology was computed, including knots with up
to 12 and links with up to 11 crossings, have 3--torsion, with the exception of
the Hopf links, whose homology is torsion-free. This is reminiscent of the omnipresence
of 2--torsion in Khovanov homology remarked by Shumakovitch \cite{torsion}.
Exemplary calculations also show the existence of 2--, 4--, 5-- and 8--torsion, which are rather scarce.
Small knots have torsion-free reduced homology, but large enough knots like the
$(8,5)$--torus knot have not.
\myparagraph{Reduced and unreduced homology:}
Reduced $\mathfrak{sl}_2$--homology
was conjectured
by Khovanov \cite{patterns}
 to be
determined by its unreduced counterpart; more precisely, that the rank of
reduced homology were one less than the rank of unreduced homology.
This is true for small knots, but 
Shumakovitch produced 15--crossing counterexamples \cite{shum-kho-and-app}.
For $\mathfrak{sl}_3$--homology even considering only knots with crossing
number six is enough to see that there is no linear relationship
between the ranks of reduced and unreduced homology.
\myparagraph{Previous computations:}
\foamho{} calculations are in agreement with the results of Carqueville and Murfet \cite{carqueville},
who compute reduced and unreduced rational $\mathfrak{sl}_3$--homology of all prime knots and links with
up to six crossings.

\bibliography{main}
\end{document}